\documentclass{article}

\makeatletter
\title{Oversampling errors in multimodal medical imaging are due to the Gibbs effect}

\author{Davide Poggiali $^{1}$, Diego Cecchin $^{1,2}$,\\
 Cristina Campi$^3$ and Stefano De Marchi $^{1,4}$\\
 \mbox{ }\\
  { \footnotesize $^{1}$ \quad Padova Neuroscience Center, University of Padova (Italy)}\\
  { \footnotesize $^{2}$ \quad Nuclear Medicine Unit, Department of Medicine - DIMED, Padua University Hospital, Italy}\\
  { \footnotesize $^{3}$ \quad Department of Mathematics, University of Genova, Italy}\\
  { \footnotesize $^{4}$ \quad Department of Mathematics ``Tullio Levi-Civita'', University of Padova, Padova (Italy)}
 }

\usepackage{amsmath}
\usepackage{amsfonts}
\usepackage{amssymb}
\usepackage{amsthm}
\usepackage[latin1, utf8]{inputenc}
\usepackage[T1]{fontenc}
\usepackage[english]{babel}
\usepackage{setspace}

\newtheorem{theorem}{Theorem}
\newtheorem{lemma}{Lemma}
\newtheorem{corollary}{Corollary}

\usepackage[colorlinks=false]{hyperref}
\hypersetup{
colorlinks  = true,
linkcolor   = black,
anchorcolor = black,
citecolor   = black,
filecolor   = black,
urlcolor    = black,
}

\usepackage{setspace}
\usepackage{pifont}
\usepackage{xcolor}
\newcommand{\red}[1]{\textcolor{red}{#1}}

\usepackage{siunitx}
\usepackage{booktabs}

\usepackage{graphicx}

\usepackage{bm}
\usepackage[normalem]{ulem}
\RequirePackage{natbib}
\begin{document}
\maketitle
\abstract{To analyse multimodal 3-dimensional medical images, interpolation is required for resampling which - unavoidably - introduces an interpolation error. In this work we describe the interpolation method used for imaging and neuroimaging and we characterize the Gibbs effect occurring when using such methods. In the experimental section we consider three segmented 3-dimensional images resampled with three different neuroimaging software tools for comparing undersampling and oversampling strategies and to identify where the oversampling error lies.
The experimental results indicate that undersampling to the lowest image size is advantageous in terms of mean value per segment errors and that the oversampling error is larger where the gradient is steeper, showing a Gibbs effect.}

\section{Introduction}

In the context of multimodal medical imaging~\cite{Ehman2017, Marti-Bonmati2010, Zhang2017}, image data of the same physical body are obtained from different imaging systems. The resulting images have different geometrical resolutions in terms of Full Width at Half Maximum (FWHM)~\cite{Cecchin2015} and by consequence with different samplings of the same Field Of View (FOV). In many cases, two types of image sources are involved, producing:
\begin{itemize}
    \item one (or more) {\bf morphological}, high-resolution image, usually obtained by Computed Tomography (CT) or by Magnetic resonance Imaging (MRI);
    \item one (or more) {\bf functional}, low-resolution image. A functional image is typically obtained by Single Photon Emission Computed Tomography (SPECT), Positron Emission Tomography (PET)~\cite{Zhang2017}, functional MRI (fMRI)~\cite{Millman2007} or by emerging systems as the Magnetic Particle Imaging (MPI)~\cite{DeMarchi2017}.
\end{itemize}

Thanks to their high spatial resolution, morphological images can be used for identifying different structures of the physical body under exam through segmentation. In order to estimate the mean activity (or any other statistical moment) of the functional images inside each of the segments or Volumes of Interest (VOIs), it is mandatory that the segmentation and the functional images have the same size. This can be achieved either by oversampling the functional image to reach the same resolution of the morphological reference or by undersampling the segmentation image~\cite{Tustison2014, Dumitrescu2019}. Despite what common sense may suggest, the latter is preferable, due to the bigger interpolation errors occurring in oversampling. This represents a paradox (we will later refer to such effect as {\it the resampling paradox}), but also a big waste of time since an accurate segmentation image at high resolution (in the case of the human brain) is obtained typically after 10-20 hours of manual work or 3-10 hours of automatic segmentation pipelines~\cite{Delgado2014}.\\

In this paper we describe the interpolation methods used in imaging and medical imaging, focusing on the Gibbs effect~\cite{Jerri1998, Lehmann1999, Fischl2002, DeMarchi2017}. The Gibbs effect results in an interpolation error due to a discontinuity and consisting in non-physical oscillations of higher amplitude around the discontinuity. We found a pointwise error bound for the proposed interpolation methods that does not vanish in the proximity of the point of discontinuity as the number of samples increases and that is larger in the intervals of sampling points closer to the discontinuity.\\

We then consider three of the most used software suites in neuroimaging performing undersampling and oversampling of three different test images. The results confirm the resampling paradox and point out that the oversampling interpolation errors are due to the Gibbs effect as the voxelwise error is mostly concentrated around the borders of the VOIs, where most of the discontinuities lie.\\
The authors of this paper refer to their experience in the field of PET/MRI neuroimaging~\cite{Cecchin2017}, but the analysis and conclusions are applicable for any multimodal image setting and any scanned physical body.\\

The paper is organized as follows. The next section is dedicated to mathematical formulations and preliminary definitions, in which we mention some results that will be useful to understand the experimental settings and the error measures. Section three is dedicated to interpolation methods in imaging along with error estimations related to Gibbs effect. In the ``Materials and Methods'' section we introduce and shortly describe the three test images used in experiments, the tools used to perform undersampling and oversampling, the error measures defined for comparing undersampling versus oversampling and last the measures defined for detecting the oversampling error spatial location. Results and some explanatory comments are presented in Sections 4 and 5. Possible overcomings to the interpolation errors issues analyzed in this paper are considered in Section 6 as well as with some recent findings which have shown their efficiency in preventing the Gibbs effect.\\

\section{Preliminaries and definitions}

\subsection{The Gibbs effect}\label{Gibbs}

The Gibbs effect is the non-physical oscillation generated when a discontinuous function is approximated by a truncated Fourier series~\cite{Jerri1998,Fornberg2011}. Such effect appears also in interpolation and does not vanish as the number of samples increases, but the undershoots and overshoots tend to remain stable. The overshoots have higher intensity around the borders of the discontinuities and tends to vanish as the border gets more distant. \\
The Gibbs phenomenon appears everywhere in signal processing; in MRI it is also referred as ``ringing effect'' (or ringing artifact)~\cite{Lehmann1999, Fischl2002, chloa20} for the wave-like oscillation appearing radially around the borders of the image discontinuities that be observed in image reconstruction and image resampling.\\

\begin{figure}[h!]
\begin{center}
   \includegraphics[width=.9\linewidth]{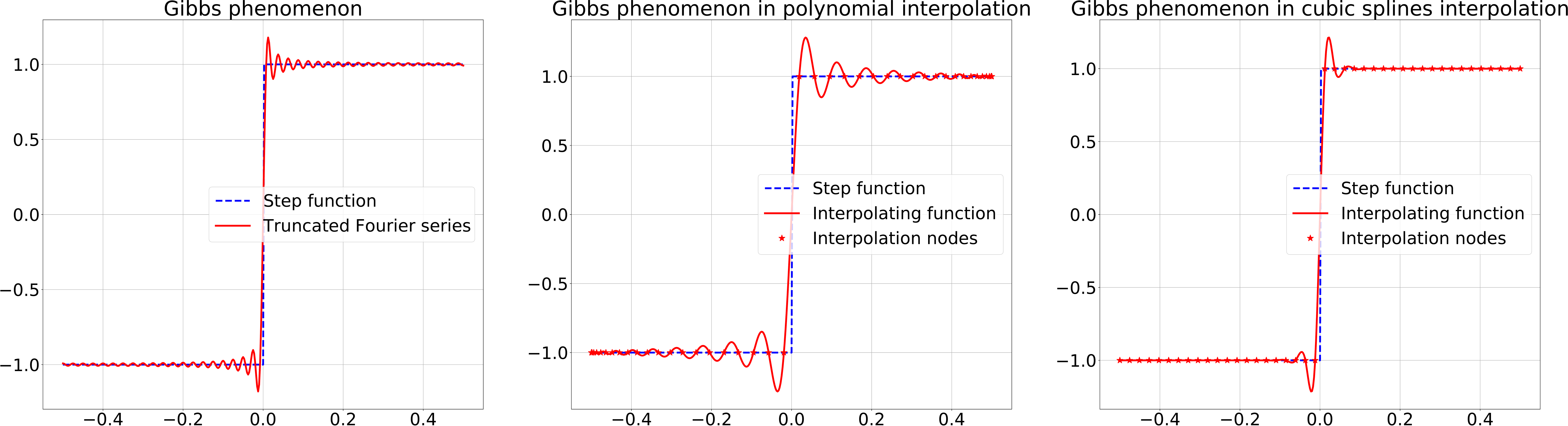}
   \caption{From left to right: Truncated Fourier series, polynomial interpolation (on Chebyshev-Lobatto nodes) and cubic splines interpolation (on equispaced nodes) of the Heaviside step function.\label{fig0}}
\end{center}
\end{figure}

\begin{figure}[h!]
\begin{center}
   \includegraphics[width=.45\linewidth]{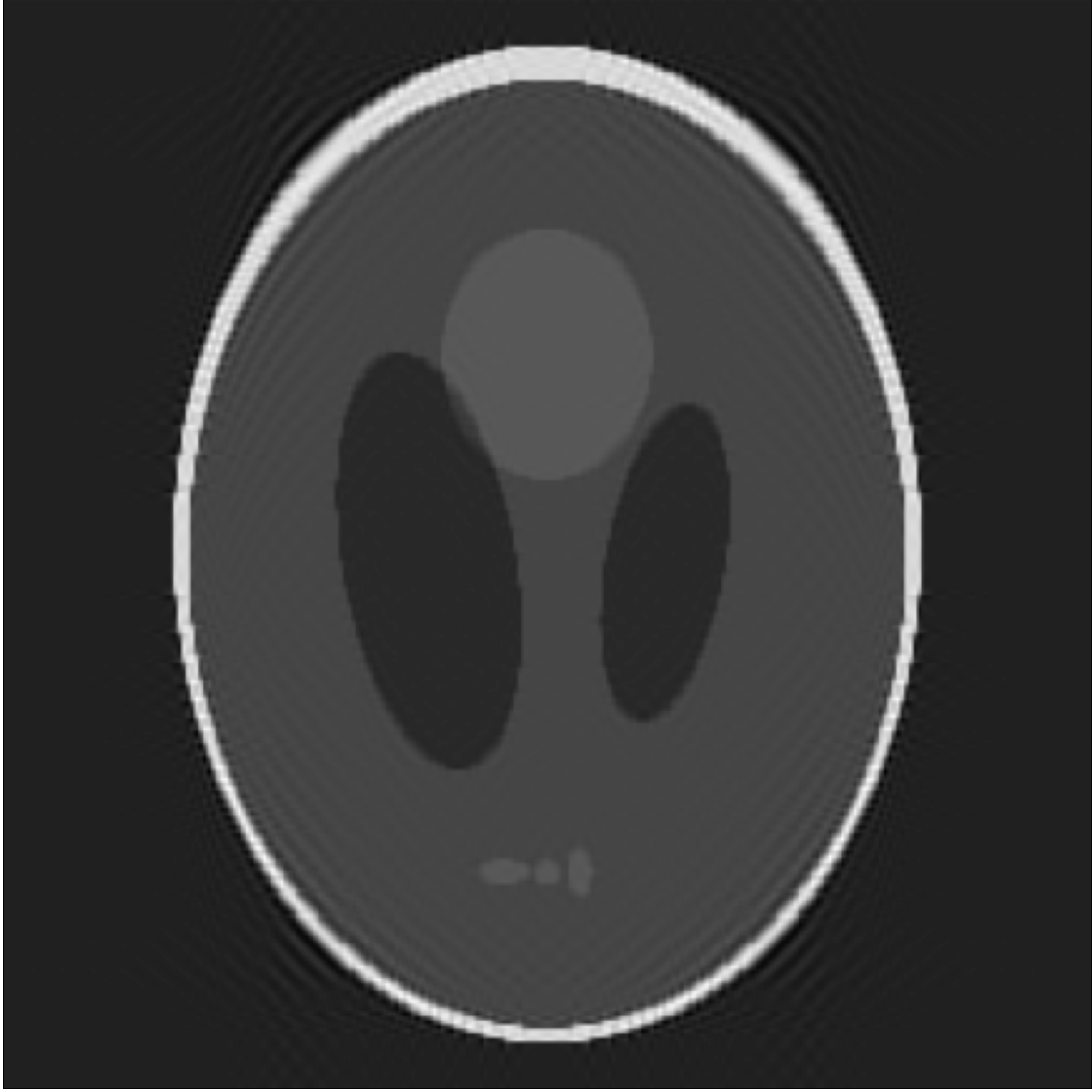}
   \caption{Ringing effect in a Shepp-Logan phantom. The effect has been magnified for visualization purposes.\label{fig0.5}}
\end{center}
\end{figure}

In Figure~\ref{fig0} it is possible to observe a simple example of Gibbs effect on the Heaviside step function; both in Truncated Fourier approximation and interpolation, oscillations not present in the original signal appear around the discontinuity, gradually mitigating their amplitude as the distance from the discontinuity increases. In Fig.~\ref{fig0.5} an example of ringing effect is given using a phantom image.

\subsection{An image and its sampling}
A 3-dimensional {\bf image} can be defined as a trivariate function~\cite{M2AN_1992__26_1_51_0}
\begin{equation}
    f:\Omega\subset\mathbb{R}^3 \to \mathbb{R}
\label{image}
\end{equation}

where the domain $\Omega$ is called Field of View (FOV) and is a rectangular parallelepiped $\Omega = [a_1,b_1] \times [a_2,b_2] \times [a_3,b_3]$.

A {\bf raster image} is obtained by sampling an image $f$ on $\Omega$ with a regular and equispaced 3-dimensional grid $X=\{{\bm x}_{ijk}\}$ of cardinality $I\times J\times K$, i.e.
\begin{equation}
    {\bm x}_{ijk} = \left(x_i, y_j, z_k\right) = \left((b_1-a_1)\frac{i-1}{I-1} + a_1, (b_2-a_2)\frac{j-1}{J-1} + a_2, (b_3-a_3)\frac{k-1}{K-1} + a_3\right)
    \label{equigrid}
\end{equation}
with $i = 1, \ldots, I, \;j = 1, \ldots, J, \;k = 1, \ldots, K.$

The {\bf intensity values} of the raster image are stored in an array $F = f|_{X}$ of size $I\times J\times K$
\begin{equation}
    F_{ijk} = f({\bm x}_{ijk}).
\label{simage}
\end{equation}
The sampling grid is usually not stored in a file to save disk space and kept implicit.\\

\subsection{Resampling an image}
{\bf Resampling} an image means to compute it over another regular and equispaced grid $\tilde{X}$ (called evaluation grid). We refer to {\bf undersampling} if $\tilde{X}$ has a lower cardinality than $X$, {\bf oversampling} in the opposite case. Being the exact function $f$ not known in general, to perform resampling it is necessary to compute the interpolant $\mathcal{P}f$ of the image intensity values on the known grid $X$ over the evaluation grid $\tilde{X}$. This means, once chosen a basis $\mathcal{B} = \{B_{ijk}(\cdot)\}$, to define an interpolating function
\begin{equation}
  \mathcal{P}f:\Omega\subset\mathbb{R}^3 \to \mathbb{R}
\end{equation}
\begin{equation}
  \mathcal{P}f(\cdot) = \sum_{ijk} c_{ijk} B_{ijk}(\cdot)
\end{equation}

with $c_{ijk}$ coefficients such that the interpolation conditions

\begin{equation}
   \mathcal{P}f (\bm{x}_{ijk}) = f(\bm{x}_{ijk})
 \label{interp_cond}
\end{equation}
are fulfilled.
At last, the resampled image is computed as the array $\mathcal{P}f|_{\tilde{X}}$.\\
Interpolation methods used in neuroimaging will be discussed in detail in the next section.

\subsection{Boolean images and morphological operators}
A {\bf Boolean image} is an image sampled over a grid as in Eq.~\eqref{simage}, but with all Boolean entries, i.e.
\begin{equation}
    M\in \{0,1\}^{I\times J\times K};
\end{equation}
a Boolean image can be represented equivalently by the set of its true-valued indices

\begin{equation}
    \mathcal{I}(M) := \{ (i,j,k) \,|\, M_{ijk} = 1 \} \subset \mathbb{Z}^3.
\label{idx}
\end{equation}

This representation allows to make some definitions (cf e.g. ~\cite{Haralick}). Given $M\in \{0,1\}^{I\times J\times K}$ a Boolean image, $I = \mathcal{I}(M)$ its index representation, and and image $F\in \mathbb{R}^{I\times J \times K}$ sampled over the same grid, we can define:

\begin{itemize}
    \item the {\bf volume} of $M$, as the cardinality of the corresponding index representation
    \begin{equation}
        \text{Vol}(M) = |\mathcal{I}(M)|
    \label{vol}
    \end{equation}
    \item the {\bf restriction of} $F$ to $M$ as the array of the values of $F$ in the voxels where $M$ is equal to 1
    \begin{equation}
        F[M] := \{F_{ijk}\, | \, (i,j,k)\in \mathcal{I}(M) \}
    \label{restr}
    \end{equation}
    \item the {\bf shift} by a vector $v\in\mathbb{Z}^3$
     $$I_v = \{ a+v \,|\, a\in I  \}$$
     \item the {\bf dilation} by a structuring element $H\subset \mathbb{Z}^3$
     \begin{equation}
       I \oplus H = \bigcup_{h\in H} I_h
     \label{dilate}
     \end{equation}
     \item the {\bf erosion} by a structuring element $M\subset \mathbb{Z}^3$
     \begin{equation}
         I \ominus H = \bigcap_{h\in H} I_{-h}
     \label{erode}
     \end{equation}
\end{itemize}

The structuring element used in this work is the 3-dimensional cross
$$H = \{ (0,0,0), (1,0, 0), (-1,0, 0), (0,1,0), (0,-1,0), (0,0,1), (0,0,-1)\}.$$

Intuitively, the dilation of a Boolean image $A$ is a Boolean image returning 1 if any of the surrounding voxel of $A$ is 1, 0 otherwise; on the other side the erosion returns 1 if all the surrounding voxels are 1, 0 otherwise.

\subsection{Segmentation of an image}
A {\bf segmentation} of an image means partitioning the FOV $\Omega$ into non-overlapping sets (cf e.g. ~\cite{Pham2000}) called {\bf VOIs} (Volumes Of Interest) from now on

\begin{equation}
    \Omega = \bigcup_{p=0}^{n} \Gamma_p\;\;\;\;\text{ such that }\; \Gamma_p\cap \Gamma_q = \emptyset \;\;\forall\;\; p\neq q = 0,\ldots,n.
\label{segmentation}
\end{equation}

For each $\Gamma_p$ we can define a Boolean image $M_p$ as the array of elements
$$(M_p)_{ijk} = \left\lbrace\begin{array}{cc}
    1 & {\bm x}_{ijk} \in \Gamma_p \\
    0 & \text{otherwise}
\end{array}\right.$$
with size equal to the cardinality of its sampling grid $X$.

The {\bf segmentation image} is an array of the same size as the image $F$ with as value the integer number $p$ corresponding to the number of segment $\Gamma_p$ to whom each voxel belongs that is

$$ M := \sum_{p=0}^n p\, M_p.$$

The zero-indexed VOI usually represents the background, the part of the FOV that does not contain the object under exam.

\subsection{Statistical moments of an image inside a VOI}

Computing a {\bf statistical moment} $\mu_N$ of an image $f$ on a VOI $\Gamma_p$ corresponds to computing the integral
$$\mu_N(f; \Gamma_p) := \frac{1}{\text{meas}(\Gamma_p)}\int_{\Gamma_p} {\bm x}^N f({\bm x}) \;d{\bm x} $$
being $N = (n_1, n_2, n_3)$ a multi-index~\cite{HuMIng, Klesk}.
Since the images are discrete, the moment is approximated by its discrete version (with a slight abuse of notation, we call it $\mu_N$ as well)
\begin{equation}
  \mu_N (f; \Gamma_p) \approx \mu_N(F[M_p]) := \frac{1}{\text{Vol}(M_p)}\sum_{(i,j,k) \in \mathcal{I}(M_p)} {\bm x}_{ijk}^N F_{ijk}
 \label{mean}
\end{equation}

being $F[M_p]$ the restriction of $F$ to $M_p$ defined in Eq.~\eqref{restr}, $\text{Vol}(M_p)$ the volume cf. Eq.\eqref{vol}, $\mathcal{I}(M_p)$ the index representation as defined in Eq.~\eqref{idx} and ${\bm x}_{ijk}$ belonging to the sampling grid $X$. cf. Eq~\eqref{simage}.

This operation is feasible only if the image and the segmentation image are sampled at the same grid $X$. To get the moment of two differently-sampled images, it is necessary to resample either the image $F$ or the VOI image $M$ to the same size of the other one.\\
Since the mean $\mu_{\bm 0} = \mu_{(0,0,0)}$ is by far the most used statistical moment in neuroimaging, we choose it for the experiments in this paper.

\section{Image interpolation and Gibbs effect}
As stated in subsection 2.2, image resampling is computed by interpolation. This section is dedicated to the definition of the most used interpolation methods in imaging and to prove some theorems that will allow us to characterize the Gibbs effect when using such methods.

\subsection{Basis construction}

Given an image
$$f:\Omega = [a_1,b_1] \times [a_2,b_2] \times [a_3,b_3] \to \mathbb{R}$$
sampled over an equispaced grid $X = \{\bm{x}_{ijk}\}$ as in Eq.~\eqref{equigrid}, we compute the interpolant as a linear combination of the elements of a chosen basis $\mathcal{B}$.

Since the size of an image can be very large (nowadays up to $10^8$ voxels), it is necessary to build the interpolation basis as cardinal, compactly supported and separable~\cite{burger_burge_2009, Moraes2020, Getreuer2011}. The basis constructions is as follows.

\begin{itemize}
    \item[(I)] Let $\omega_a: \mathbb{R}_+ \to \mathbb{R}$, with $a\in\mathbb{N}$ the {\bf support radius}, a function such that:
\begin{itemize}
   \item[(i)] $\omega_a(0) = 1$;
   \item[(ii)] $\omega_a(k) = 0\;\; \forall~ k \in \mathbb{N}\setminus \{0\}$;
   \item[(iii)] $\omega_a(r) = 0\;\;\forall~ r \geq a$.
\end{itemize}
A list functions satisfying (i)-(iii) is given in Table~\ref{tab0} (more functions  can be found in~\cite{burger_burge_2009, Getreuer2011}).

\item[(II)] Let $t_i = (\beta - \alpha)\frac{i}{N} + \alpha$, with $i = 0, \ldots, N$ a set of equispaced points of the interval $[\alpha, \beta]$. Then we define for $t\in {[\alpha,\beta]}$ the univariate basis
\begin{equation}
   w_{\sigma}(t - t_i) = \frac{\omega_a(\sigma|t-t_i|)}{\sum_{j=0}^N \omega_a(\sigma |t-t_j|)}.
\label{basis}
\end{equation}
with $\sigma = \frac{N}{\beta-\alpha}$. \\
This basis is {\bf cardinal} for (i)-(ii), has {\bf compact support} for (iii), {\bf radial} as it depends from the absolute value of its argument, and {\bf normalized} for construction.

\item[(III)] The trivariate interpolation basis corresponding to the set $X$ is
$$\mathcal{B} = \left\lbrace  W_{ijk}(\cdot) \right\rbrace_{ijk} $$
where $W_{ijk}$ is the {\bf separable} kernel function
$$W_{ijk}(x, y, z) = w_{\sigma_1}(x-x_i)\; w_{\sigma_2}(y-y_i) \;w_{\sigma_3}(z-z_i)$$
with $\sigma_1 = \frac{b_1-a_1}{I-1}$, $\sigma_2 = \frac{b_2-a_2}{J-1}$, and $\sigma_3 = \frac{b_3-a_3}{K-1}$.
\end{itemize}

It is easy to verify that the basis $\mathcal{B}$ is {\bf cardinal} over $X$, i.e.
\begin{equation}
    W_{ijk}(\bm{x}_{i'j'k'} ) = \left \lbrace\begin{array}{cc}
   1  & \text{ if } i'=i, j'=j, k'=k \\
   0  &  \text{ otherwise}
\end{array}\right.
\label{cardinal}
\end{equation}

Furthermore, this basis is normalized by Eq.~\eqref{basis}

\begin{equation}
   \sum_{ijk} W_{ijk}(\bm{x}) = 1\;\;\;\forall \bm{x} \in \Omega.
\label{norm3}
\end{equation}

\subsection{Interpolation by convolution}

Under the assumptions of the previous subsection, we define {\bf interpolant} at $\bm{x} = (x,y,z)\in\Omega$
\begin{equation}
\mathcal{P}_f (\bm{x}) = \sum_{ijk} f(\bm{x}_{ijk}) W_{ijk}(\bm{x}).
\label{interp3}
\end{equation}

Since the basis is cardinal, the interpolation conditions~\eqref{interp_cond} are fulfilled
$$\mathcal{P}_f (\bm{x}_{i'j'k'}) = \sum_{ijk} f(\bm{x}_{ijk}) W_{ijk}(\bm{x}_{i'j'k'}) = f(\bm{x}_{i'j'k'})\;\;\;\forall i',j',k'.$$

Being $W$ a separable kernel function, Eq.~\eqref{interp3} translates to

\begin{equation}
\mathcal{P}_f (\bm{x}) = \sum_{k=1}^K \left( \sum_{j=1}^J  \left(\sum_{i=1}^I f(\bm{x}_{ijk}) w_{\sigma_1}(x-x_i)\right) w_{\sigma_2}(y-y_j)\right) w_{\sigma_3}(z-z_k)
\end{equation}
 and the evaluation of the interpolant at a given point can be performed in three steps:
 \begin{enumerate}
     \item compute $\alpha_{jk} = \sum_{i=1}^I f(\bm{x}_{ijk}) w_{\sigma_1}(x-x_i)$;
     \item then get $\beta_k = \sum_{j=1}^J  \alpha_{jk} w_{\sigma_2}(y-y_i)$;
     \item at last $\mathcal{P}_f (\bm{x}) = \sum_{k=1}^K \beta_k w_{\sigma_3}(z-z_i)$.
 \end{enumerate}

Since the basis function are separable, the interpolant~\eqref{interp3} is equivalent to a discrete convolution and the literature names it {\it interpolation by convolution}~\cite{burger_burge_2009}.

It is interesting to notice that such a technique is easily extendable to larger dimensions. Since each dimension is treated separately, we will discuss the interpolation problem as univariate. The results concerning error estimates that we present in the next subsection are easy to extend to the trivariate interpolation by convolution introduced above.

 \subsection{Error estimates for the univariate interpolation}
  Let us suppose without loss of generality that the function
$$f:[0,1] \to \mathbb{R}$$
is known at $x_i = \frac{i}{N},\;\;i=0,\ldots,N$.

By using a function $\omega_a(r)$ satisfying (I), we define the interpolating function

\begin{equation}
    \mathcal{P}^N_f(x) = \sum_{i=0}^N w(x-x_i) f(x_i).
    \label{interp1D}
\end{equation}
with $w$ as in~\eqref{basis} with $\sigma$ set to $N$.

The following shows the advantages of using compact support functions to save memory and CPU time.
\begin{lemma}
Suppose that $x\in (x_k, x_{k+1})$ for some $k\in 0, \ldots, N-1$ and that the support radius is $a\in\mathbb{N}$. Then
$$\mathcal{P}^N_f(x) = \sum_{i=\max(0, k-a)}^{\min(N, k+a-1)} w(x-x_i) f(x_i).$$

\label{lemma1}
\end{lemma}

\begin{proof}
If $i<k-a$ it will result
$$N|x-x_i| = N( (x_k -x_i) + (x-x_i)) = N\left( \frac{k-i}{N} + (x-x_i)\right) > a.$$ In similar way it is easy to show that if  $i> k+a-1$ then $N|x-x_i|>a$.\\
By consequence of (iii) in both cases $\omega_a(N|x-x_i|)=0$ and hence $w(x-x_i) = 0$.
\end{proof}

This means that - at most - only the $a$ nodes before and the $a$ nodes after $x$ are involved in computing $\mathcal{P}^N_f(x)$.\\

In order to get an estimate of the error we now define two constants that will be useful for overestimating the terms $|w(x-x_i)|$. Let

\begin{equation}
    M_{w} = \max_{x\in [0,1]} \left| \omega_a(N|x-x_i|)\right|,
    \label{M1}
\end{equation}
and
\begin{equation}
    m_{w} = \min_{x\in [0,1]} \left|\sum_{j=0}^N \omega_a(N|x-x_j|) \right|.
    \label{m1}
\end{equation}

Both constants are relatively easy to compute analytically. In fact, supposing $x\in(x_k,x_{k+1})$ and $a < k < N-a$, we can exploit the radiality of $w$ and the fact that the nodes are evenly spaced to rewrite the sum in Eq.~\eqref{m1} by coupling the $a$ nodes before and after the point $x$
\begin{equation}
    \sum_{j=k-a}^{k+a-1} \omega_a(N|x-x_j|)  = \sum_{p=0}^{a-1} \omega_a(\delta_k + p) + \omega_a(1-\delta_k + p),
    \label{m2}
\end{equation}
where $\delta_k = N|x-x_k|\in (0,1)$. Hence, we can compute exactly the constant as
\begin{equation}
    m_{w} = \min_{t\in [0,1]} \left|\sum_{p=0}^{a-1} \omega_a(t + p) + \omega_a(1-t + p) \right|.
    \label{m3}
\end{equation}

More easily, Eq.~\eqref{M1} can be rewritten as
\begin{equation}
    M_{w} = \max_{t\in [0,a]} \left| \omega_a(t)\right|.
    \label{M2}
\end{equation}

\begin{table}[h!]
\begin{center}
  \caption{Name and definition of the basis functions $\omega_a$ used in imaging and neuroimaging, along with the support radius $a$ and the constants $m_{\omega}$ and $M_{\omega}$. The shape parameter for the Gaussian in this paper is chosen as $\varepsilon=2$; in such case $m_w \approx 0.736$. As $\text{sinc}$ function we mean the normalized sinc, defined as $\text{sinc}(x) = sin(\pi x) / (\pi x)$ if $x\neq 0$; $\text{sinc}(0)=0$.\label{tab0}}

  \begin{tabular}{c l | c c c}
  Name & Definition $\omega_a(r)$ & $a$ & $m_{w}$ & $M_{w}$\\
  \toprule
  (Truncated) Gaussian & $\begin{cases} \text{exp}(-\varepsilon^2 r^2) & \text{ if } r\in [0,1)\\0&\text{otherwise}\end{cases}$ & 1& $2\text{exp}(-\varepsilon^2/4)$ & 1\\ \midrule
   Nearest Neighbour & $\begin{cases} 1 & \text{ if } r\in [0,1/2)\\0&\text{otherwise}\end{cases}$ & 1& 1& 1\\ \midrule
   Linear  & $\begin{cases} 1-r & \text{ if } r\in [0,1)\\0&\text{otherwise}\end{cases}$ & 1& 1& 1\\ \midrule
  Cubic  & $\begin{cases} r^3 -2r^2 +1 & \text{ if } r\in [0,1)\\-r^3 +5r^2 -8r+4 & \text{ if } r\in [1,2)\\0&\text{otherwise}\end{cases}$ & 2& 0.75& 1\\ \midrule
  Lanczos & $\begin{cases} \text{sinc}(x) \text{sinc}(x/2) & \text{ if } r\in [0,2)\\0&\text{otherwise}\end{cases}$ & 2& $1+2\text{sinc}(3/2)\text{sinc}(3/4)\approx 0.8726$& 1\\
  \bottomrule
  \end{tabular}
\end{center}
\end{table}

Once calculated the values of $m_{w}$ and $M_{w}$ we can bound
\begin{equation}
     |w(x-x_i)| \leq \frac{M_{w}}{m_{w}}
    \label{basis_bound}
\end{equation}
for every $x\in (x_k, x_{k+1})$ such that $a\leq k \leq N-a+1 $ and $i =k-a,\ldots,k+a-1$.\\

The values of $m_{w}$ and $M_{w}$ have been computed exactly for each choice of $\omega_a$ and shown in Table~\ref{tab0}. The reader, if interested can see the computation of the constants in the GitHub page of this paper \url{https://github.com/pog87/GibbsEffectMultimodal}. It is worth noticing that in general it is not true that $m_{w}\neq 0$; such condition will be assumed as hypothesis in what follows.\\

We can finally enounce and prove a theorem for the pointwise interpolation error bound. The theorem can be seen as a particular case of the results published in~\cite{Wendland}, of which we give a different proof. Other error estimates for Gibbs effect in radial basis function interpolation can be found in~\cite{Fornberg2011, Jung2007}.

\begin{theorem}[Pointwise error bound]
Let $f:[0,1]\to \mathbb{R}$ be a function sampled over $N+1$ equispaced points $x_i = \frac{i}{N}, \;i=0,\ldots,N$, and $w$ the basis function defined in Eq.~\eqref{interp1D}. If the interpolation point $x\in(x_k, x_{k+1})$ with $a\leq k \leq N-a+1 $ and $m_{w}>0$, then
\begin{equation}
|\mathcal{P}^N_f (x) - f(x)| \leq \frac{M_{w}}{m_{w}} \sum_{i=k-a}^{k+a-1} |f(x) - f(x_i)|.
    \label{ineq}
\end{equation}
\label{thm1}
\end{theorem}

\begin{proof}
Being the basis normalized to one and by effect of lemma~\ref{lemma1}
$$\mathcal{P}^N_f (x) - f(x) = \sum_{i=k-a}^{k+a-1} w(x-x_i)(f(x_i) - f(x)).$$
Taking the absolute values, by triangular inequality and Eq.~\eqref{basis_bound} we get~\eqref{ineq}.
\end{proof}

This result ensures that for $\frac{M_{w}}{m_{w}}$ small enough the error will not be propagated uncontrollably. Looking at Table~\ref{tab0}, we can see that such constant ratio is sufficiently small for each of the chosen basis functions. A particular attention must be dedicated in choosing the shape parameter for the Gaussian-based interpolation. In fact if we choose $\varepsilon$ too large the $m_{w}$ consequently goes to zero. In the experiments of this paper the shape parameter is set to $\varepsilon = 2$ in accordance to the most popular neuroimaging software.\\




\begin{corollary}
In the hypothesis of Theorem~\ref{thm1}, if $f$ is continuous in $x\in[0,1]$, then
$$\lim_{N\to\infty} |\mathcal{P}^N_f (x) - f(x)| = 0 .$$
\label{corol1}
\end{corollary}
\begin{proof}
If $x=0$ or $x=1$ the result holds. Otherwise we can assume without loss of generality that for $N$ big enough, $x\in(x_k, x_{k+1})$ is not close to the domain boundaries, i.e. $a \leq k \leq N-a$.\\
Therefore
$$ \sum_{i=k-a}^{k+a-1}|f(x) - f(x_i))| =  \sum_{i=k-a}^{k}|f(x) - f(x_i))| + \sum_{j=k+1}^{k+a-1}|f(x) - f(x_j))|$$
passing to the limit we get
\begin{equation}
    a\left|f(x) - f(x^-)\right| + a\left|f(x) - f(x^+)\right|
   \label{ErrLim}
\end{equation}
Being $f$ continuous in $x$, $f(x^-) = f(x^+) = f(x)$, hence by Theorem~\ref{thm1} we conclude.
\end{proof}

As consequence of this corollary, the interpolation scheme described in this work will not be subject to Runge effect. This holds for any local interpolation scheme as e.g. splines.\\

\begin{corollary}
In the hypothesis of Theorem~\ref{thm1}, if $f$ is discontinuous at $x\in[0,1]$, and the limit of the error exists, then
$$0 \leq \lim_{N\to\infty} |\mathcal{P}^N_f (x) - f(x)| \leq a\frac{M_{w}}{m_{w}} |f(x^-) - f(x^+)| .$$
\label{coroll1.5}
\end{corollary}

This shows that in case of discontinuity, the error bound does not vanish when increasing the number of samples, which corresponds to the description of the Gibbs effect given in subsection~\ref{Gibbs}.\\

As last theoretical result, we show that the error bound increases as the evaluation point $x$ get closer to the discontinuity. Let us suppose that
\begin{equation}
   f(t) = \begin{cases} f^-(t) & t\leq \xi \\ f^+(t) & t> \xi\end{cases}
   \label{disc}
\end{equation}
with $f^-\in \mathcal{C}(\xi - \epsilon, \xi)$ and $f^+\in \mathcal{C}(\xi, \xi + \epsilon)$ for some $\epsilon>0$ and that the evaluation point $x$ belongs to $(\xi-\epsilon, \xi + \epsilon)$. As showed before, if $N|x-\xi|>a$, the discontinuity has no effect on $\mathcal{P}^N_f$, as the interpolation acts locally. The next Corollary concerns what happens in the opposite case.

\begin{corollary}
In the hypothesis of Theorem~\ref{thm1}, if $f$ is bounded and discontinuous at $\xi$, $x\in (x_k, x_{k+1})$ with $a < k < N-a$ and $\xi \in [x_{k+p}, x_{k+p+1}]$ with $p \in\{ -a, \ldots, a-2\}$, calling
$$ F = \sum_{k-a}^{k+a-1}|f(x)- f(x_i)|,$$
then:\\
\begin{itemize}
    \item {\bf Case~1} \hspace*{.5cm} $\xi\leq x$ (so $p\leq 0$):
\begin{equation}
    pc_1 + c_2 \leq F \leq pC_1 + C_2
    \label{less}
\end{equation}
with $c_1,c_2, C_1, C_2\in \mathbb{R}$, $c_1, C_1>0$; \vspace*{.7cm}

    \item{\bf Case~2} \hspace*{.5cm} $\xi> x$ (so $p>0$):
    \begin{equation}
    pr_1 + r_2 \leq F \leq pR_1 + R_2
    \label{more}
\end{equation}
with $r_1,r_2, R_1, R_2\in \mathbb{R}$, $r_1, R_1<0$.
\end{itemize}
\label{corol2}
\end{corollary}

\begin{proof}
Let us suppose that $f^-((\xi - \epsilon, \xi)) = [a,b]$ and $f^+((\xi, \xi+\epsilon)) = [\alpha, \beta]$ with $a\leq b<\alpha\leq \beta$. This assumption will not make loose generality, in fact if $N$ is large enough the images of the left and right neighbourhoods will not intersect, and in the opposite case $\alpha\leq\beta<a\leq b$ the proof is similar. Now, let
$$\Delta f^+ = \beta-\alpha;\;\;\;\Delta f^- = b-a;$$
$$\Delta = \beta-a;\;\;\;\delta = \alpha-b.$$
Let us split $F$ in the $(a+p+1)$ nodes before and the $(a-p-1)$ nodes after the discontinuity, so that
$$F = \sum_{i=k-a}^{k+p}|f(x)- f(x_i)| + \sum_{j=k+p+1}^{k+a-1}|f(x)- f(x_j)|.$$
In {\bf Case~1} the $x_j$ will be on the same side of $x$ and the $x_i$ on the opposite, so all the terms $|f(x) - f(x_j)|$ can be underestimated with zero and overestimated with $\Delta f^+$, while all the terms $|f(x)- f(x_i)|$ will fall between $\delta$ and $\Delta$; hence
$$ (a+p+1) \delta \leq F \leq (a+p+1)\Delta + (a-p-1)\Delta f^+ = p ( \Delta- \Delta f^+) + (a+1)\Delta + (a-1)\Delta f^+ ,$$
proving the first case as $\Delta > \Delta f^+$ and $\delta>0$.\\

Similarly for the {\bf Case~2},
$$ (a-p+1) \delta \leq F \leq (a+p+1)\Delta f^+ + (a-p-1)\Delta,$$
which proves~\eqref{more}.
\end{proof}

By Corollary~\ref{corol2}, the value of $F$ will increase as $p$ approaches zero both from the left or from the right, and by consequence the error bound shown in Theorem~\ref{thm1} will increase with the decrease of $|p|$, where $|p|$ is the distance between the intervals of nodes in which $x$ and $\xi$ lie. This corresponds to the description of the Gibbs Effect in subsection~\ref{Gibbs}, stating that the overshoots and undershoots can be larger in proximity of a discontinuity.

\subsection{Error estimates for the trivariate interpolation}

The error estimate given in the previous subsection can be extended to the trivariate case - which corresponds to an image interpolation by convolution - and in general to any dimensional separable interpolation.\\

\begin{theorem}
Let $f$ be a 3d image~\eqref{image} sampled over an equispaced grid $X$~\eqref{equigrid}.
Let $\mathcal{P}_f$ any interpolant of form~\eqref{interp3} using basis functions $$W_{ijk}(x, y, z) = w_{\sigma_1}(x-x_i)\; w_{\sigma_2}(y-y_i) \;w_{\sigma_3}(z-z_i)$$ built using (I)-(III) and with $w_{\sigma_1}, w_{\sigma_2}, w_{\sigma_3}$ such that $m_w>0$, then
\begin{equation}
 |\mathcal{P}_f(\bm{x}) - f(\bm{x})| \leq \left( \frac{M_w}{m_w}\right)^3 \sum_{ijk}|f(\bm{x}) - f(\bm{x}_{ijk})|
 \label{thm2}
\end{equation}
\end{theorem}

\begin{proof}
Since the basis is normalized~\eqref{norm3}, we can write
$$\mathcal{P}_f(\bm{x}) - f(\bm{x}) = \sum_{ijk} \left( f(\bm{x}_{ijk}) - f(\bm{x}) \right) W_{ijk}(\bm{x}).$$
By taking the absolute value, using the triangular inequality and applying Theorem~\ref{thm1} three times, we conclude.

\end{proof}

The interested reader can see the experimental verification of this error bound in the GitHub repo \url{https://github.com/pog87/GibbsEffectMultimodal}.

\section{Materials and Methods}

\subsection{Images}
A set of three images have been used for numerical experiments

\begin{enumerate}
    \item The 3D Shepp-Logan (SL) phantom~\cite{Shepp1974}, a picewise-constant function made by the weighted sum of characteristic functions over different ellipsoids. The SL images have been created in python with {\tt tomopy} and {\tt nibabel}~\cite{Brett2020} and saved in nifti format. The Segmention was made with a python script by grouping the voxels of the same intinsity values. Size: $256\times 256\times 256$;
    \item A 3D isotropic MRI, T1-weighted of one of the authors' head. Skull-stripping and automatical segmentation was made with {\tt GIF}~\cite{Cardoso2015a}. Size: $180\times 560\times 560$;
    \item The CT of a walnut, downloaded from\\
    \href{http://www.informatik.uni-leipzig.de/~wiebel/public_data/index.html}{http://www.informatik.uni-leipzig.de/$\sim{}$wiebel/public\_data/}\\
    along with its segmentation image~\cite{Prassni2010}. Size: $400\times 296\times 352$.
\end{enumerate}

\begin{figure}[h!]
\begin{center}
   \includegraphics[width=.85\linewidth]{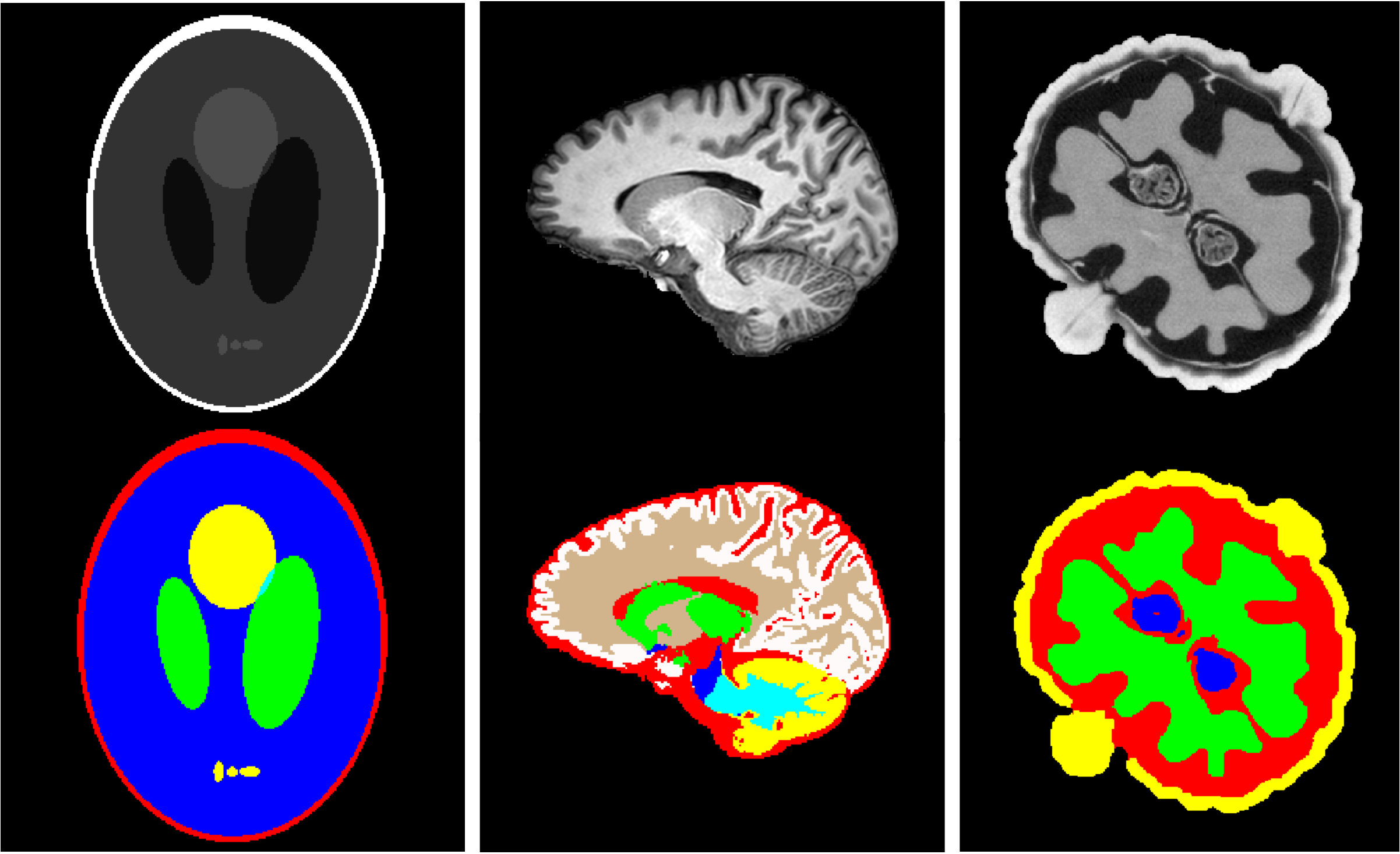}
   \caption{From left to right: a slice of the the SL phantom, Brain MRI and walnut CT. In the second row the corresponding segmentation images, with a different color used for each different VOIs.\label{fig1}}
\end{center}
\end{figure}

In Fig.~\ref{fig1} we observe a slice of each image and the corresponding segmentation image.

Each image has been undersampled with a factor 2 per dimension using {\tt nibabel}, except for the SL phantom; in this case the exact function of the image is known and hence the undersampled version has been calculated analytically. This corresponds to the common PET/MRI setting, where the morphological MRI image is about $(1\;mm)^3$ per voxel whereas the functional PET image is sampled at $(2\;mm)^3$.\\

The undersampled version of each image in this experiment acts as a functional image. We decide not to take into account the Partial Volume Effect (PVE)~\cite{Muller1992} of the imaging systems because the aim of this work is to investigate only the interpolation errors in resampling. Since the high-resolution values are known, we are able to estimate the interpolation errors that occurs in resampling.

\subsection{Software used}
Three different tools from commonly used neuroimaging software suites have been used to perform interpolation for undersampling and oversampling.
\begin{enumerate}
    \item {\tt antsApplyTransforms} from Advanced Normalization Tools ({\bf ANTs})~\cite{AVANTS2008, Avants2011} v2.2.0, with interpolation options:
    \begin{enumerate}
        \item {\tt Nearest} and {\tt Multilabel} for undersampling segmentation;
        \item {\tt Linear}, {\tt Gaussian}, {\tt LanczosWindowedSinc} and {\tt Splines} for oversampling {\it functional} image;
    \end{enumerate}
    \item {\tt flirt} from FMRIB Software Library ({\bf FSL})~\cite{Jenkinson2012,Reuter2010, Smith2004} v5.0.8, with interpolation options:
    \begin{enumerate}
        \item {\tt Nearest} for undersampling;
        \item {\tt Trilinear} and {\tt Splines} for oversampling;
    \end{enumerate}
    \item {\tt mri\_convert} from {\bf Freesurfer}~\cite{Fischl2002, Zou2004} v5.3, with interpolation options:
    \begin{enumerate}
        \item {\tt Nearest} for undersampling;
        \item {\tt Trilinear} and {\tt Splines} for oversampling.
    \end{enumerate}
\end{enumerate}

The {\tt Nearest}, {\tt (Tri-)Linear}, and {\tt Splines} options refer to the Nearest Neighbour, Linear and Cubic basis function respectively, as reported in Table~\ref{tab0}. {\tt Multilabel} means that each VOI has been Gaussian-filtered, interpolated with cubic splines and finally each voxel has been assigned to the argument of maximal value. {\tt Gaussian} and {\tt LanczosWindowedSinc} indicate a Gaussian and Lanczos-2 interpolation, respectively. The formulation of each of these functions can be found in Table~\ref{tab0}.

The whole dataset and code used for the experiments has been uploaded to the page \url{https://github.com/pog87/GibbsEffectMultimodal} for the sake of reproducibility.

\subsection{Evaluating the error}

As indicator of the performance of the undersampling versus oversampling procedures we choose the relative 2-norm error with respect to the reference mean values of the original, high-resolution image.\\
The reference mean values are defined as the mean intensity values of the original high-resolution image $I$ and the high-resolution segmentation $M$,
    $$v = \mu_{\bm 0}(I[M]) = \{\mu_{\bm 0}(I[M_p])\}_{p=0,\ldots,n},$$
    being $\mu_{\bm 0}$ the discrete mean~\eqref{mean}. Such values are available in our experiments because of the choice to use undersampled images acting as functional, but are unknown in real cases of multimodal image analysis .

We computed the undersampling and oversampling approximations of $v$ as the means:
\begin{enumerate}
    \item by using functional image $F$ and undersampled segmentation $M_{-}$
    $$v_{-} = \mu_{\bm 0}(F[M_{-}]);$$
    \item by using oversampled functional image $F_{+}$ and segmentation M
    $$v_{+} = \mu_{\bm 0}(F_{+}[M]);$$
\end{enumerate}

Hence, the relative 2-norm error between the reference values and the under/over-sampling values is evaluated for each image and interpolation method as

\begin{equation}
    \text{err}_{\pm} = \frac{ \| v - v_{\pm} \|_2}{\| v  \|_2}
\label{rel_err}
\end{equation}

Undersampling and oversampling errors will be compared for all the test images and software to confirm or deny the resampling paradox. Results are presented in the dedicated section.\\

\subsection{Locating the oversampling error}

\begin{figure}[h!]
\begin{center}

\includegraphics[width=.8\linewidth]{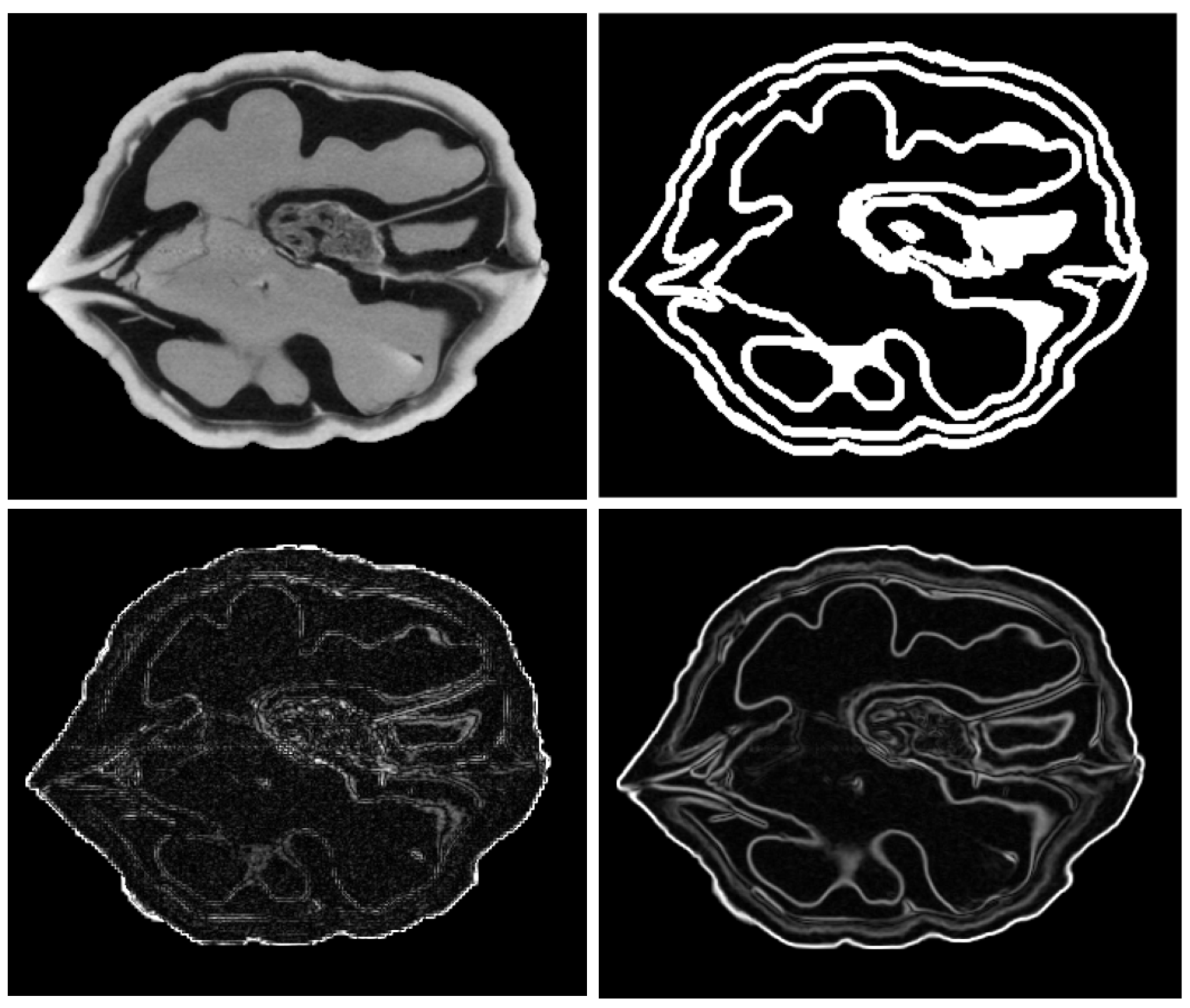}
\caption{Upper left: a slice of the walnut CT $I$. Upper right: the borders $\delta M$ computed as in Eq.~\eqref{borders} at the same slice. Bottom left: the pointwise absolute error $E$ (as in Eq.~\eqref{voxelwise_err}) in ANTs linear interpolation. Bottom right: $\nabla I$, the norm of the gradient of the original image computed with Eq.~\eqref{nabla}. \label{fig2}}
\end{center}
\end{figure}

As the segmentation of high-resolution images is time-consuming, it is preferable to use the oversampling procedure, but it could be a source of larger errors. To better understand such interpolation error, it is necessary to further investigate the source of the oversampling interpolation error. Namely, we want to know where the voxelwise error

\begin{equation}
E = |I - F_{+}|
\label{voxelwise_err}
\end{equation}

is larger. As an example, in Fig.~\ref{fig2} (bottom left) a slice of the voxelwise interpolation error is shown.\\

We computed a VOI around the borders of the segmentation image as

\begin{equation}
 \delta M = \vee_{p=0}^n \delta M_p
\label{borders}
\end{equation}
where $\vee$ means the logic {\tt OR} operator, and each border $\delta M_p$ is computed as
$$\delta M_p = \text{\tt dilate}(M_p,3) \otimes \text{\tt erode}(M_p,3)$$
where $\otimes$ denotes the voxelwise logic {\tt xor} operator, ${\tt dilate}(X,n)$ and ${\tt erode}(X,n)$ the binary erosion an dilation morphological operators repeated $n$ times. Erosion and dilation operators are implemented in {\tt scikit-image} Python package~\cite{VanDerWalt2014}. An example of $\delta M$ can be found in Fig.~\ref{fig2} (top right).\\

We also estimate the 2-norm of the gradient of the high-resolution image $I$

\begin{equation}
    \nabla I := \left\| \left( \frac{\partial f}{\partial x}, \frac{\partial f}{\partial y}, \frac{\partial f}{\partial z} \right) \right\|_2
\label{nabla}
\end{equation}

where each partial derivative is estimated with {\it Sobel operator} (cf. ~\cite{Kittler}) in {\tt scikit-image}. An example of $\nabla I$ is in Fig.~\ref{fig2} (bottom right).

As a dissimilarity measure between images, we use the Dissimilarity Structural SIMilarity index ({\bf DSSIM}), defined as
$$\text{DSSIM}(A,B) := 1- \text{SSIM}(A,B)$$
where the Structural SIMilarity index (or {\bf SSIM}, is a $[0,1]$-valued similarity measure which has been shown to be a stable metric for multivariate interpolation~\cite{wangSSIM, Dumitrescu2019, Marchetti2021} and it is defined as
$$ \text{SSIM}(A,B) = \frac{(2\bar{A}\bar{B} + c_1)(2\Sigma_{AB} + c_2)}{(\bar{A}^2 + \bar{B}^2 + c_1)(\Sigma_A^2 + \Sigma_B^2 + c_2)}$$
being $\bar{A}$ the mean, $\Sigma_A$ the standard deviation, and $c_1,c_2$ two real default constants used for stability. For the tuning of the parameters $c_1, c_2$ refer e.g. to~\cite{wangSSIM}. The SSIM of two identical images is 1, which means that the DSSIM is zero.\\

Using {\tt scikit-image} we compute the global error $\text{DSSIM}(I,  F_{+})$,  the error at the borders $\delta M$ as
$\text{DSSIM}(I[\delta M],  F_{+}[\delta M])$ and the percent ratio between the two in order to assess how much of the error is located in the borders. \\
Moreover, we want to know how much of the image gradient is located at the borders and which part of the FOV volume is taken by the border in terms of number of voxels.

In fact although we expect an ideal segmentation to perfectly identify all of the different parts of the object in the FOV, in practice a segmentation can miss some signal discontinuity. We decided to take into account
$$\nabla I[\delta M]_{(\%)} := \frac{\| \nabla I[\delta M] \|_1}{\|\nabla I\|_1} \cdot 100$$

as a measure of the ratio of the image gradient located in the VOIs border $\delta M$. Moreover, to ensure that $\delta M$ does not cover too much of the FOV - in such case a high gradient ratio would be unavoidable - we computed the Volume ratio

$$\text{Vol}(\delta M)_{(\%)} := \frac{\text{Vol}(\delta M)}{\text{meas}(\Omega)} \cdot 100.$$


\section{Experimental results}

The results of the VOI-wise analysis are shown in Table~\ref{tab1}. In all the cases but one, undersampling the segmentation images leads to a smaller relative 2-norm error on estimating the mean intensity in the VOIs.

\begin{table}[h!]
  \caption{Relative errors in 2-norm of the mean value per VOI, computed with Eq~\eqref{rel_err}. Minimal error per software per image in red. Empty cells indicate the same as above.\label{tab1}}

  \begin{tabular}{| c c c | S[table-auto-round, table-format=1.5]  S[table-auto-round, table-format=1.5]  S[table-auto-round, table-format=1.5] |}
  \toprule
  Software & sampling &method & {SL phantom} & {Brain MRI} & {Walnut CT} \\
  \toprule
  ANTs & under & Nearest    & 0.11715 & \red{0.00251} & \red{0.00029}\\
       &       & MultiLabel & \red{0.11714} & 0.00287 & 0.00098\\ \cline{2-6}
       & over  & Linear     & 0.16559 & 0.02427 & 0.01965\\
       &       & Gaussian   & 0.24197 & 0.05616 & 0.04605\\
       &       & Lanczos    & 0.12989 & 0.01415 & 0.01425\\
       &       & Splines    & 0.12790 & 0.00791 & 0.00751\\
  \midrule
  \midrule
  FSL & under & Nearest  & 0.10845 & \red{0.00327} & \red{0.00383}\\ \cline{2-6}
      & over & Trilinear & 0.16528 & 0.02430 & 0.02025\\
      &      & Splines   & \red{0.10403} & 0.00792 & 0.00759\\
  \midrule
  \midrule
  Freesurfer & under & Nearest  & \red{0.11716} & \red{0.00252} & \red{0.00029}\\ \cline{2-6}
             & over & Trilinear & 0.16559 & 0.02427 & 0.01965\\
             &      & Splines   & 0.12644 & 0.00792 & 0.00751\\
  \bottomrule
  \end{tabular}
\end{table}

In Table~\ref{tab2} the DDSIM and DSSIM at the borders are listed for each image, software and interpolation method. The ratio of the error located in the borders is  80 to 85\% for the SL phantom, whose borders host the totality of the gradient, being the phantom a picewise-constant function. The brain image holds the 89\% of the gradient and an error percentage of about 30-55\%. The walnut image has the lowest gradient percentage (71\%) and the lowest error percentage, spanning from 19 to 31\%.

\begin{table}[h!]
\begin{center}
  \caption{Interpolation oversampling errors in terms of $\text{DSSIM}(I,  F_{+})$ (the smaller the better), DSSIM at the segmentation border VOI $\text{DSSIM}(\delta M) = \text{DSSIM}(I[\delta M],  F_{+}[\delta M])$ and its percentage. Per image, per method. Percentage of the image gradient at the borders $\delta M$ and percentage volume of $\delta M$ over the whole FOV. Empty cells indicate the same as above.\label{tab2}}

  \begin{tabular}{lll | S[table-auto-round, table-format=1.3] S[table-auto-round, table-format=1.3]  S[table-auto-round, table-format=3.1]|S[table-auto-round, table-format=3.1]S[table-auto-round, table-format=3.1]}
  \toprule
     Image & Software &   method & {DSSIM} &  {DSSIM($\delta M$)} &  {\%} & {$\nabla I[\delta M]_{(\%)}$} & {$\text{Vol}(\delta M)_{(\%)}$} \\
  \midrule
       SL &     ANTs &  Linear &  0.037153 &              0.031298 &   84.242829 &     100.000000 &    14.897966 \\
          &          &    Gaussian &  0.047003 &          0.038652 &   82.233787 &                                &          \\
          &          &     Lanczos &  0.042190 &          0.035423 &   83.961040 &                                &          \\
          &          &     Splines &  0.041255 &          0.034414 &   83.418321 &                                &          \\  \cline{2-8}
          &      FSL &  Trilinear &  0.029491 &          0.023908 &   81.068121 &                                &          \\
          &          &     Splines &  0.033593 &          0.027387 &   81.524186 &                                &          \\ \cline{2-8}
          &    Freesurfer &  Trilinear &  0.037153 &          0.031298 &   84.242829 &                                &          \\
          &          &     Splines &  0.038853 &          0.031272 &   80.488007 &                                &          \\ \midrule

               Brain&     ANTs &  Linear &  0.062901 &          0.034254 &   54.457968 &      89.20926 &    18.886478 \\
          &          &    Gaussian &  0.111262 &          0.051932 &   46.675099 &            &          \\
          &          &     Lanczos &  0.080985 &          0.033880 &   41.835420 &            &          \\
          &          &     Splines &  0.102826 &          0.033518 &   32.596547 &            &          \\ \cline{2-8}
          &      FSL &  Trilinear &  0.062913 &          0.034261 &   54.457173 &            &          \\
          &          &     Splines &  0.102840 &          0.033526 &   32.599943 &            &          \\ \cline{2-8}
          &    Freesurfer &  Trilinear &  0.062901 &          0.034254 &   54.457961 &            &          \\
          &          &     Splines &  0.102826 &          0.033518 &   32.596537 &            &          \\ \midrule

   Walnut &     ANTs &  Linear &  0.044889 &          0.011284 &   25.137363 &      71.41374 &    14.234087 \\
          &          &    Gaussian &  0.091228 &          0.028299 &   31.020122 &        &     \\
          &          &     Lanczos &  0.079741 &          0.021932 &   27.503712 &        &     \\
          &          &     Splines &  0.114494 &          0.021974 &   19.191816 &        &     \\ \cline{2-8}
          &      FSL &  Trilinear &  0.046841 &          0.012816 &   27.360825 &        &     \\
          &          &     Splines &  0.114795 &          0.022474 &   19.577444 &        &     \\ \cline{2-8}
          &    Freesurfer &  Trilinear &  0.044889 &          0.011284 &   25.137374 &        &     \\
          &          &     Splines &  0.114494 &          0.021973 &   19.191790 &        &     \\
  \bottomrule

  \end{tabular}
\end{center}
\end{table}

\section{Discussion}

The resampling paradox is confirmed by this data shown in Table~\ref{tab1}, as the errors in estimating the mean value of the VOIs assessed by undersampling the segmentation image is smaller than any error produced by oversampling the {\it functional} images. The only exception is the splines interpolation in FSL, which seems to be produced by an error compensation in the computation of the mean values. This is confirmed by its DSSIM shown in Table~\ref{tab2}, which is larger than the DSSIM given by the Trilinear interpolation in FSL.\\

This data also indicates that the interpolation error occurring in oversampling is a Gibbs effect, confirming the results proved in the theoretical section. As we can intuitively infer from Fig~\ref{fig2}, the pointwise error is mostly located in voxels where the gradient is higher. This intuition is confirmed by inspecting Table~\ref{tab2} as the most of the error lies at the borders of the VOIs when the gradient at the border is higher and decreases in cases where the segmentation misses some high-gradient zones.\\
In the case of the walnut CT image, it can be observed in Fig.\ref{fig1} and in Fig.\ref{fig2} that the provided segmentation misses some of the inner skin surrounding the seed and some gaps within the seed.\\

It is also noticeable that the error percentage in the borders has the tendency to be lower if the error is higher, possibly showing a propagation of the Gibbs effect around the borders.\\

\section{Conclusions and future works}

The Gibbs effect has been proved to appear in interpolation by convolution and is confirmed by experimental results to be the largest source of interpolation error in multimodal medical imaging.\\

The resampling paradox is confirmed experimentally. Undersampling is - as a matter of fact - the most chosen option in multimodal neuroimaging, despite the fact that the segmentation at high resolution results in a waste of effort and time.\\
In order to avoid such effect a new interpolation technique is needed allowing an oversampling of functional images which minimizes the Gibbs effect. A promising approach is given by the interpolation by convolution with the scale factor Point Spread Function (sfPSF)~\cite{Cardoso2015} which takes into account the different FWHM of the morphological and functional images. In alternative, a spectral filtering~\cite{DeMarchi2017a} could be a fast and efficient way of dealing with the Gibbs effect.
Another interesting method we can consider is the {Fake-Nodes} interpolation introduced in \cite{DeMarchi2020}, which has shown to be an effective approach for univariate interpolation without resampling and multivariate approximation of data ~\cite{DeMarchi2020a, DeMarchi2021}, reducing the Gibbs effect.

This research has been funded by the PNC - Padova Neuroscience Center, University of Padova (Italy) as part of the project "A computational tool for neurodegenerative stratification using PET/RM" and partially funded by GNCS-IN$\delta$AM.




This research has been accomplished within Rete ITaliana di Approssimazione (RITA).

The authors have no conflicts of interest to declare.


\bibliography{main.bbl}

\end{document}